\documentclass[12pt,a4paper,oneside,reqno]{article}

\usepackage[utf8]{inputenc}
\usepackage[english]{babel}
\usepackage[left=1in,right=1in,top=1in,bottom=1in]{geometry}
\usepackage{amsfonts,amssymb,amsmath,amsthm}
\usepackage{enumerate}
\usepackage{cite}
\usepackage{xcolor}
\usepackage[normalem]{ulem}
\usepackage{longtable}
\usepackage{hyperref}
\hypersetup{colorlinks=true,linkcolor=blue,citecolor=blue,urlcolor=blue}
\usepackage{fancyhdr}
\theoremstyle{plain}
\newtheorem*{ThmOne}{Theorem 1 \cite{Kac77}}
\newtheorem*{ThmTwo}{Theorem 2 \cite{Kac77}}
\newtheorem*{LemA}{Lemma A}
\newtheorem*{LemB}{Lemma B}
\theoremstyle{remark}
\newtheorem*{RemC}{Remark C (on Theorem 2(c))}
\newtheorem*{RemD}{Remark D (on the use of Lemma 6)}
\newtheorem*{RemE}{Remark E (the case $B(1,1)$)}

\newcommand{\FF}{F}
\newcommand{\Antider}{\mathrm{Antider}}
\newcommand{\spn}[1]{\operatorname{span}\langle #1\rangle}
\newcommand{\eps}{\epsilon}

\newcommand{\doi}[1]{\href{https://doi.org/#1}{doi:#1}}

\begin{document}

\noindent\fbox{\parbox{\dimexpr\textwidth-2\fboxsep-2\fboxrule\relax}{\small
\textbf{Corrected version verified in Lean.} This is a corrected version of the paper
I.~B.~Kaygorodov, \emph{$\delta$-Derivations of classical Lie superalgebras}, Siberian Mathematical
Journal 50 (2009), no.~3, 434--449, \doi{10.1007/s11202-009-0049-x}.}}

\bigskip

\begin{center}
{\Large $\delta$-derivations of classical Lie superalgebras}

\bigskip

{\bf Ivan Kaygorodov}
\end{center}

\bigskip

\noindent {\bf Abstract.} {\it We consider the $\delta$-derivations of classical Lie superalgebras and
prove that these superalgebras admit nonzero $\delta$-derivations only when $\delta=0,\frac12,1$. The
structure of $\frac12$-derivations for classical Lie superalgebras is completely determined.}

\medskip


\bigskip


\bigskip

\setcounter{tocdepth}{1}
\tableofcontents

\bigskip

\section*{Introduction}
\addcontentsline{toc}{section}{Introduction}

Hopkins \cite{Hop96} considered the antiderivations of Lie algebras. An antiderivation is a special
case of a $\delta$-derivation --- a linear mapping $\mu$ of an algebra such that
$\mu(xy)=\delta(\mu(x)y+x\mu(y))$, where $\delta$ is some fixed element of the ground field.
Independently, these results were obtained in a generalized form by Filippov \cite{Fil98}. He proved
that a prime Lie $\Phi$-algebra equipped with a nondegenerate symmetric invariant bilinear form does not
admit nonzero $\delta$-derivations if $\delta\notin \{-1,0,\frac12,1\}$. In the same article, some description
of $\frac12$-derivations was given for an arbitrary prime Lie $\Phi$-algebra $A$ ($\frac16\in\Phi$)
equipped with a nondegenerate symmetric invariant bilinear form. He proved that a linear mapping
$\varphi\colon A\to A$ is a $\frac12$-derivation if and only if $\varphi\in\Gamma(A)$, where $\Gamma(A)$
is the centroid of $A$. Let $A$ be a simple Lie algebra over a field of characteristic $p\neq2,3$.
Assume that $A$ is equipped with a nondegenerate symmetric invariant bilinear form. It follows from the
results cited above that every $\frac12$-derivation $\varphi$ of $A$ is of the shape
$\varphi(x)=\alpha x$, $\alpha\in\Phi$. Later Filippov  described the $\delta$-derivations for
prime alternative and non-Lie Malcev $\Phi$-algebras under some restrictions on the ring of operators
$\Phi$. He proved that the algebras in these classes do not admit nonzero $\delta$-derivations if
$\delta\notin\{0,\frac12,1\}$ \cite{Fil00}.

In \cite{Kay07}, the $\delta$-derivations were studied for simple finite-dimensional Jordan superalgebras
over an algebraically closed field of characteristic $0$ and for semisimple finite-dimensional Jordan
algebras over an algebraically closed field of characteristic distinct from $2$. It was proved that these
classes of algebras and superalgebras possess nonzero non-trivial $\delta$-derivations only when
$\delta=\frac12$. A complete description of $\frac12$-derivations was given for these
classes, and it was shown that $\varphi$ is a $\frac12$-derivation if and only if $\varphi(x)=\alpha x$
for some $\alpha\in\FF$ (for the simple superalgebras and the simple algebras).

In the present work we give a description of nontrivial $\delta$-derivations of classical Lie
superalgebras. We prove that the superalgebras in this class admit nonzero $\delta$-derivations only if
$\delta=0,\frac12,1$. We give a complete description of $\frac12$-derivations for this class of
superalgebras. We prove that $\varphi$ is a $\frac12$-derivation if and only if $\varphi(x)=\alpha x$ for
some $\alpha\in\FF$.

\section{Basic facts and definitions}

We denote by $\spn{a,b}$ the linear span of $a$ and $b$. Let $U$ be a vector space, and let $V$ be a
subspace in $U$. Take $x\in U$ and denote by $x|_V$ the projection of $x$ to $V$ (with respect to a
fixed complement of $V$, which is always clear from the context).

Let $\FF$ be an algebraically closed field of characteristic $0$. A Lie superalgebra $G=G_0+G_1$ is a
$\mathbb Z_2$-graded algebra with the superidentities
\begin{equation}\label{eq1}
[x,y]=-(-1)^{p(x)p(y)}[y,x],\qquad [[x,y],z]-[x,[y,z]]-(-1)^{p(y)p(z)}[[x,z],y]=0,
\end{equation}
where $p(x)=i$ if $x\in G_i$. The second superidentity is a generalization of the Jacobi identity.

A finite-dimensional Lie superalgebra $G=G_0+G_1$ is called classical provided that $G$ is simple and
the representation of $G_0$ on $G_1$ is completely reducible. In \cite{Kac77}, Kac gave a complete
classification of classical Lie superalgebras over an algebraically closed field of characteristic $0$.

\begin{ThmOne}
Let $A$ be a classical Lie superalgebra over an algebraically closed field of characteristic $0$. Then
$A$ is one of the superalgebras
\[
A(m,n),\ A(n,n),\ B(m,n),\ D(m,n),\ C(n),\ P(n),\ Q(n),\ D(2,1;\beta),\ G(3),\ F(4).
\]
\end{ThmOne}

We now recall a description of the superalgebras mentioned in Theorem~1.

$A(m,n)$: Let $sl(m,n)$ be the following subspace in the matrix superalgebra $M_{m+n,m+n}$ with the
$\mathbb Z_2$-grading:
\begin{longtable}{lcl}
$(sl(m,n))_0 $&$=$&$\left\{\begin{pmatrix}A&0\\0&D\end{pmatrix} : A\in M_m(\FF),\ D\in M_n(\FF),\
\mathrm{tr}(A)=\mathrm{tr}(D)\right\},$\\ 
$(sl(m,n))_1$&$=$&$\left\{\begin{pmatrix}0&B\\C&0\end{pmatrix} : B\in M_{m,n}(\FF),\ C\in M_{n,m}(\FF)\right\}.$
\end{longtable}
Equip $sl(m,n)$ with the structure of a Lie superalgebra by $[a,b]=ab-(-1)^{p(a)p(b)}ba$. If $m=n$ then
this superalgebra contains the one-dimensional ideal $\langle E_{2n}\rangle$ consisting of the scalar
matrices $\lambda E_{2n}$, $\lambda\in\FF$. The Lie superalgebra $sl(1,1)$ is three-dimensional and
nilpotent. We put 
\begin{longtable}{lcl}
$A(m,n)$&$=$&$sl(m+1,n+1)$ for $m\neq n$, $m,n\geq0$,\\ $A(n,n)$&$=$&$sl(n+1,n+1)/\langle E_{2n+2}\rangle$ for $n>0$.
\end{longtable}

$B(m,n)$, $D(m,n)$, $C(n)$, $P(n)$, and $Q(n)$ are some subsuperalgebras in $A(k,l)$ (for $Q(n)$: a
quotient of a subsuperalgebra by its centre); recall that $(B(m,n))_0=B_m\oplus C_n$,
$(D(m,n))_0=D_m\oplus C_n$ ($m\ge2$), $(C(n))_0=\FF\oplus C_{n-1}$, and that the $G_0$-module $G_1$ is
$V\otimes W$ (the tensor product of the natural modules) for $B(m,n)$, $C(n)$, $D(m,n)$,
$V\otimes W^*\oplus V^*\otimes W$ for $A(m,n)$, $S^2V\oplus\Lambda^2V^*$ for $P(n)$
($(P(n))_0=A_n$) and the adjoint module for $Q(n)$ ($(Q(n))_0=A_n$) \cite{Kac77}.

$D(2,1;\alpha)$, $\alpha\in\FF^*\setminus\{0,-1\}$: This is a one-parameter family of $17$-dimensional
Lie superalgebras consisting of all simple Lie superalgebras such that $(D(2,1;\alpha))_0$ is a Lie algebra
of type $G_1^1\oplus G_1^2\oplus G_1^3$ ($G_1^j\cong A_1$) and its representation on $(D(2,1;\alpha))_1$
is $sl_2\otimes sl_2\otimes sl_2$.

$F(4)$: Define $F(4)$ as a $40$-dimensional classical Lie superalgebra such that $(F(4))_0$ is a Lie
algebra of type $B_3\oplus A_1$ and its representation on $(F(4))_1$ is $spin_7\otimes sl_2$.

$G(3)$: Define $G(3)$ as a $31$-dimensional classical Lie superalgebra for which $(G(3))_0$ is a Lie
algebra of type $G_2\oplus A_1$ and its representation on $(G(3))_1$ is $G_2\otimes sl_2$.

The classical superalgebras distinct from $Q(n)$, $P(n)$, and $A(1,1)$ were called the basic classical
Lie superalgebras in \cite{Kac77}. We will use this terminology.

Below we describe the systems of roots for the basic classical Lie superalgebras \cite{Kac77}. Denote by
$\Delta_0$ and $\Delta_1$ the systems of even and odd roots respectively. Denote a simple system of roots
by $\Pi$. In these cases a Cartan subalgebra $H$ is a subspace of the space $D$ of diagonal matrices. The
roots are expressed in terms of the standard basis $\eps_i$ on the dual space $D^*$ (more accurately, the
restrictions of $\eps_i$ on $H$).

$A(m,n)$: The root system [in terms of  $\eps_1,\dots,\eps_{m+1}$,
$\delta_1=\eps_{m+2},\dots,\delta_{n+1}=\eps_{n+m+2}$]:
\begin{longtable}{lcl}
$\Delta_0$&$=$&$\{\eps_i-\eps_j;\ \delta_i-\delta_j\},\ i\neq j;$\\ 
$\Delta_1$&$=$&$\{\pm(\eps_i-\delta_j)\};$\\
$\Pi$&$=$&$\{\eps_1-\eps_2,\dots,\eps_{m+1}-\delta_1,\delta_1-\delta_2,\dots,\delta_n-\delta_{n+1}\}.$
\end{longtable}

$B(m,n)$: The root system [in terms of $\eps_1,\dots,\eps_m$,
$\delta_1=\eps_{2m+1},$ $\dots,$ $\delta_n=\eps_{2m+n}$]:

\begin{longtable}{lcl}
$\Delta_0$&$=$&$\{\pm\eps_i\pm\eps_j;\ \pm2\delta_i;\ \pm\eps_i;\ \pm\delta_i\pm\delta_j\},\ i\neq j;$\\
$\Delta_1$&$=$&$\{\pm\delta_i;\ \pm\eps_i\pm\delta_j\};$\\
$\Pi$&$=$&$\{\delta_1-\delta_2,\dots,\delta_n-\eps_1,\eps_1-\eps_2,\dots,\eps_{m-1}-\eps_m,\eps_m\},\ m>0,$\\
$\Pi$&$=$&$\{\delta_1-\delta_2,\dots,\delta_{n-1}-\delta_n,\delta_n\},\ m=0.$
\end{longtable}

$C(n)$: The root system [in terms of $\eps=\eps_1$, $\delta_1=\eps_3,\dots,
\delta_{n-1}=\eps_{n+1}$]:
\begin{longtable}{lcl}
$\Delta_0$&$=$&$\{\pm2\delta_i;\ \pm\delta_i\pm\delta_j\};$\\  $\Delta_1$&$=$&$\{\pm\eps\pm\delta_i\};$\\
$\Pi$&$=$&$\{\delta_1-\delta_2,\dots,\delta_{n-2}-\delta_{n-1},\delta_{n-1}-\eps,\delta_{n-1}+\eps\}.$
\end{longtable}

$D(m,n)$: The root system [in terms of $\eps_1,\dots,\eps_m$,
$\delta_1=\eps_{2m+1},\dots,\delta_n=\eps_{2m+n}$]:
\begin{longtable}{lcl}
$\Delta_0$&$=$&$\{\pm\eps_i\pm\eps_j;\ \pm2\delta_i;\ \pm\delta_i\pm\delta_j\},\ i\neq j;$\\
$\Delta_1$&$=$&$\{\pm\eps_i\pm\delta_j\};$\\
$\Pi$&$=$&$\{\eps_1-\eps_2,\dots,\eps_m-\delta_1,\delta_1-\delta_2,\dots,\delta_{n-1}-\delta_n,2\delta_n\}.$
\end{longtable}

$D(2,1;\alpha)$: The root system [in terms of $\eps_1,\eps_2,\eps_3$]:
\begin{longtable}{lcl}
$\Delta_0$&$=$&$\{\pm2\eps_1,\pm2\eps_2,\pm2\eps_3\}$;\\
$\Delta_1$&$=$&$\{\pm\eps_1\pm\eps_2\pm\eps_3\}$;\\ 
$\Pi$&$=$&$\{\eps_1+\eps_2+\eps_3,-2\eps_1,-2\eps_2\}$.
\end{longtable}
$F(4)$: The root system [in terms of the functionals $\eps_1,\eps_2,\eps_3$, corresponding
to $B_3$, and $\delta$, corresponding to $A_1$]:
\begin{longtable}{lcl}
$\Delta_0$&$=$&$\{\pm\eps_i\pm\eps_j;\ \pm\eps_i;\ \pm\delta\},\ i\neq j;$\\
$\Delta_1$&$=$&$\{\tfrac12(\pm\eps_1\pm\eps_2\pm\eps_3\pm\delta)\};$\\
$\Pi$&$=$&$\{\tfrac12(\eps_1+\eps_2+\eps_3+\delta),-\eps_1,\eps_1-\eps_2,\eps_2-\eps_3\}.$
\end{longtable}

$G(3)$: The root system [in terms of $\eps_1,\eps_2,\eps_3$, corresponding to $G_2$,
$\eps_1+\eps_2+\eps_3=0$, and $\delta$, corresponding to $A_1$]:
\begin{longtable}{lcl}
$\Delta_0$&$=$&$\{\eps_i-\eps_j;\ \pm\eps_i;\ \pm2\delta\};$\\ $\Delta_1$&$=$&$\{\pm\eps_i\pm\delta;\ \pm\delta\};$\\
$\Pi$&$=$&$\{\delta+\eps_1,\eps_2,\eps_3-\eps_2\}.$
\end{longtable}

In what follows, $G_\beta$ denotes the root subspace corresponding to the root $\beta$, and $g_\beta$
denotes an element of this space (a nonzero one, unless stated otherwise).

For the basic classical Lie superalgebras, we have the following

\begin{ThmTwo}
Let $G$ be a basic classical Lie superalgebra, and let $G=\bigoplus G_\alpha$ be its root decomposition
with respect to a Cartan subalgebra $H$. Then

{\rm(a)} $G_0=H$,

{\rm(b)} $\dim(G_\alpha)=1$ for $\alpha\neq0$,

{\rm(c)} $[G_\alpha,G_\beta]\neq0$ if and only if $\alpha,\beta,\alpha+\beta\in\Delta$ \textup{(}for
$G\neq A(n,n)$, $n\ge2$; see Remark~C\textup{)}.
\end{ThmTwo}

Here $G_0$ in~(a) denotes the root space of the zero weight (not the even part).

\begin{RemC}
For $G=A(n,n)$, $n\ge2$, statement~(c) fails. Indeed, for $psl(3,3)$ the functionals satisfy
$\eps_1+\eps_2+\eps_3=\delta_1+\delta_2+\delta_3$ on $H$ (the supertrace vanishes on $H$, and $H$ is taken
modulo the scalar matrices), hence $(\eps_1-\delta_1)+(\eps_2-\delta_2)=\delta_3-\eps_3$ is an odd root,
while $[e_{14},e_{25}]=e_{14}e_{25}+e_{25}e_{14}=0$, where $G_{\eps_1-\delta_1}=\FF e_{14}$ and
$G_{\eps_2-\delta_2}=\FF e_{25}$. Statement~(c) is not used below.
\end{RemC}

Take $\delta\in\FF$. A linear mapping $\varphi$ of a superalgebra $A$ is called a $\delta$-derivation
provided that for arbitrary $x,y\in A$ we have
\begin{equation}\label{eq2}
\varphi(xy)=\delta(x\varphi(y)+\varphi(x)y).
\end{equation}
The definition of a $1$-derivation coincides with the conventional definition of a derivation; a
$0$-derivation is an arbitrary endomorphism $\varphi$ of $A$ such that $\varphi(A^2)=0$. By a nontrivial
$\delta$-derivation we mean a nonzero $\delta$-derivation, which is neither a $1$-derivation nor
$0$-derivation. Clearly, the multiplication operator by an element of the ground field $\FF$ is a
$\frac12$-derivation in every superalgebra. We are interested in studying the action of nontrivial
$\delta$-derivations on classical Lie superalgebras over an algebraically closed field $\FF$ of
characteristic $0$. Note that a $\delta$-derivation is an arbitrary linear map (it need not preserve
the grading), and no signs occur in~\eqref{eq2}. In a Lie superalgebra the product $xy$ is the bracket
$[x,y]$; both notations are used below.

\section{The main lemmas}

In this section we formulate and prove the main lemmas which will be used in the subsequent
investigation of the action of $\delta$-derivations of classical Lie superalgebras.

\medskip\noindent\textbf{Lemma 3.} {\it Let $\varphi$ be a  $\delta$-derivation of a Lie
superalgebra $G$ and $x\in G_1$. Then $\varphi(x^2)\in G_0$.}

\medskip\noindent\textsc{Proof.} Let $\varphi(x)=x_0+x_1$, $x_i\in G_i$. Then
\[
\varphi(x^2)=\delta([\varphi(x),x]+[x,\varphi(x)])=\delta([x_0+x_1,x]+[x,x_0+x_1])=2\delta[x_1,x]\in G_0.
\]
The lemma is proved.

\medskip\noindent\textbf{Lemma 4.} {\it Let $\varphi$ be a  $\delta$-derivation of a
classical Lie superalgebra $G=G_0+G_1$. Then $\varphi(G_0)\subseteq G_0$.}

\medskip\noindent\textsc{Proof.} In \cite{Kac77}, it was shown that $[G_1,G_1]=G_0$. Then, for an
arbitrary $x\in G_0$, $x=\sum_{i=1}^{n_x}y_iz_i$, where $y_i,z_i\in G_1$. It is easy to see that
$x=\sum_{i=1}^{n_x}\big(\frac14(y_i+z_i)^2-\frac14(y_i-z_i)^2\big)$. Now, we obtain the claim by Lemma~3.
The lemma is proved.

\medskip\noindent\textbf{Lemma 5.} {\it Let $\varphi$ be a nontrivial $\delta$-derivation of a simple Lie
algebra $G$ and $\dim(G)\geq3$ over an algebraically closed field $\FF$ of characteristic $0$. Then either $\delta=\frac12$ and
there exists $\alpha\in\FF$ such that $\varphi(x)=\alpha x$ for every $x\in G$ or $\delta=-1$ and
$G\cong A_1$.}

\medskip\noindent\textsc{Proof.} The claim follows from Theorems~3, 5, and~6 in \cite{Fil98}.

\medskip\noindent\textbf{Remark.} the hypothesis that $\FF$ is algebraically closed cannot be dropped: the
algebra $\mathfrak{sl}_2(\mathbb C)$ regarded as a $6$-dimensional Lie algebra over $\mathbb R$ is simple,
and $\varphi(x)=ix$ is a $\frac12$-derivation which is not of the form $\alpha x$ with $\alpha\in\mathbb R$.
In the paper $\FF$ is algebraically closed throughout (Section~1), so this does not affect the other
results. A complete proof of Lemma~5 (following the scheme of \cite{Fil98}: invariant form,
$\varphi\in\Gamma(G)$ for $\delta\neq-1,\frac12,1$ and so $\varphi=0$; the centroid of a central simple
algebra for $\delta=\frac12$; root decomposition and $\mathfrak{sl}_2$-strings for $\delta=-1$) is
verified in Lean~4.

\medskip

Express the Lie algebra $A_1$ in the shape of the algebra of columns of order $3$ with the product as
follows:
\[
\begin{bmatrix}a\\b\\c\end{bmatrix}\begin{bmatrix}x\\y\\z\end{bmatrix}=
\begin{bmatrix}bx-cy\\2ay-2bx\\2cx-2az\end{bmatrix},\quad\text{where}\quad
\begin{bmatrix}a&b\\c&-a\end{bmatrix}\mapsto\begin{bmatrix}a\\b\\c\end{bmatrix}.
\]
The first product is understood as $\left[\begin{bmatrix}a\\b\\c\end{bmatrix},
\begin{bmatrix}x\\y\\z\end{bmatrix}\right]=\begin{bmatrix}bz-cy\\2ay-2bx\\2cx-2az\end{bmatrix}$ (the
commutator of the corresponding matrices); i.e., $h=(1,0,0)^T$, $e=(0,1,0)^T$, $f=(0,0,1)^T$ with
$[h,e]=2e$, $[h,f]=-2f$, $[e,f]=h$. Denote by $\Antider(A_1)$ the space of antiderivations of $A_1$.

\medskip\noindent\textbf{Lemma 6.} {\it $\Antider(A_1)=\left\{\begin{bmatrix}-2a&b&c\\2c&a&d\\2b&e&a
\end{bmatrix} : a,b,c,d,e\in\FF\right\}$ \textup{(}the matrices act on the columns\textup{)}.}

\medskip\noindent\textsc{Proof.} The claim follows from \cite{Hop96}.

\medskip

In other words, a linear map $\varphi$ of $A_1=\spn{h,e,f}$ is an antiderivation if and only if
\begin{equation}\tag{$\ast$}\label{eqAnti}
\varphi(h)=-2ah+be+cf,\qquad \varphi(e)=\tfrac c2\,h+ae+df,\qquad \varphi(f)=\tfrac b2\,h+e'e+af
\end{equation}
for some $a,b,c,d,e'\in\FF$ (the first column of the matrix of Lemma~6 is $(-2a,2c',2b')^T$; put
$b=2c'$, $c=2b'$). This is the form in which Lemma~6 is used below; see Remark~D.

\begin{RemD}
In the proofs of Lemmas~12--19 of the published version Lemma~6 was used in the form
$\varphi(h)=-2ah+be+cf$, $\varphi(e)=2ch+ae+df$, $\varphi(f)=2bh+e'e+af$, i.e., with the transposed
matrix. Such maps are in general not antiderivations: for $b=1$, $a=c=d=e'=0$ we get
$\varphi([e,f])=\varphi(h)=e$, whereas $-[\varphi(e),f]-[e,\varphi(f)]=-[e,2h]=4e$. The correct
form is~\eqref{eqAnti}: the coefficients $2c$ and $2b$ must be replaced by $\frac c2$ and $\frac b2$.
All computations below are carried out with~\eqref{eqAnti}; in Lemmas~12, 13, 19 this changes
only intermediate coefficients, in Lemma~14 the case $n=2$ is reproved, and the computations of
Lemmas~15, 17, 18 are replaced by an application of Lemmas~A and~B.
\end{RemD}

\medskip\noindent\textbf{Lemma 7.} {\it Let $\varphi$ be a nontrivial $\delta$-derivation of a Lie algebra
$L=L_1\oplus L_2$, where $L_1$ is a semisimple Lie algebra \textup{(}and $\delta\neq0$\textup{)}.
Then $\varphi(L_2)\subseteq L_2$.}

\medskip\noindent\textsc{Proof.} Let $x_i\in L_i$ and $\varphi(x_i)=x_i^1+x_i^2$, where $x_i^j\in L_j$.
Then
\[
0=\varphi([x_1,x_2])=\delta\big([x_1,x_2^1+x_2^2]+[x_1^1+x_1^2,x_2]\big)=
\delta\big([x_1,x_2^1]+[x_1^2,x_2]\big),
\]
whence $[x_1,x_2^1]=0$. It is known that the center of a semisimple Lie algebra is zero. Therefore,
$x_2^1=0$ and $\varphi(L_2)\subseteq L_2$. The lemma is proved.

\medskip\noindent\textbf{Lemma 8.} {\it Let $G=G_0+G_1$ be a basic classical Lie superalgebra, and let
$\varphi$ be a nontrivial $\delta$-derivation of $G$ with $\delta\neq-\frac12$. Then
$\varphi(G_1)\subseteq G_1$.}

\medskip\noindent\textsc{Proof.} Let $H$ be a Cartan subalgebra in $G$. Take $g_\beta\in G_1\cap
G_\beta$ and $h\in H$. Then $\varphi(g_\beta)=\sum_{\gamma\in\Delta}k_\gamma^\beta g_\gamma$
$+\,h_\beta$, $h_\beta\in H$, and
\begin{multline*}
\sum_{\gamma\in\Delta}\beta(h)k_\gamma^\beta g_\gamma+\beta(h)h_\beta=\varphi(\beta(h)g_\beta)=\\
\varphi(hg_\beta)=\delta(\varphi(h)g_\beta+h\varphi(g_\beta))=
\delta\Big(\varphi(h)g_\beta+\sum_{\gamma\in\Delta}\gamma(h)k_\gamma^\beta g_\gamma\Big).
\end{multline*}
Since $\varphi(h)g_\beta\in G_1$ \textup{(}Lemma~4\textup{)} and $h\in H$ is arbitrary, we have
$h_\beta=0$ and $\beta=\delta\gamma$ when $\gamma\in\Delta_0$ and $k_\gamma^\beta\neq0$.
An odd root can be proportional to an even root $\gamma$ only if $\gamma=\pm2\beta$; hence if
$\delta\neq\frac12$ \textup{(}and $\delta\neq-\frac12$\textup{)} then the lemma is proved. If
$\delta=\frac12$ then $\varphi(g_\beta)=g_1^\beta+g_{2\beta}$ and
$\varphi(g_{-\beta})=g_1^{-\beta}+g_{-2\beta}$, where $g_1^\beta,g_1^{-\beta}\in G_1$,
$g_{2\beta},g_{-2\beta}\in G_0$ \textup{(}possibly zero\textup{)}. Then
$\frac12\big(g_{2\beta}g_{-\beta}+g_1^\beta g_{-\beta}+g_\beta g_1^{-\beta}+g_\beta g_{-2\beta}\big)=
\varphi(g_\beta g_{-\beta})\in G_0$, i.e., $g_{2\beta}g_{-\beta}+g_\beta g_{-2\beta}=0$; the two
summands have the different weights $\beta$ and $-\beta$, so both vanish, and since
$[G_{2\beta},G_{-\beta}]=G_\beta$ and $[G_\beta,G_{-2\beta}]=G_{-\beta}$ when $2\beta\in\Delta$, we get
$g_{2\beta}=g_{-2\beta}=0$, whence $\varphi(g_\beta)\in G_1$. The lemma is proved.

\medskip

The following two lemmas (new in this version) replace most of the case analysis of Section~3.

\begin{LemA}
Let $G=G_0+G_1$ be a classical Lie superalgebra, $\delta\neq0,1$, and let $\varphi$ be a
$\delta$-derivation of $G$ such that $\varphi(S)=0$ for some simple ideal $S$ of the Lie algebra $G_0$.
Then $\varphi=0$.
\end{LemA}

\noindent\textsc{Proof.} By the description of the $G_0$-module $G_1$ (Section~1 and \cite{Kac77}),
every simple ideal of $G_0$ acts on $G_1$ without nonzero invariants; since $S$ is semisimple and $G_1$
is a finite-dimensional $S$-module, $G_1=[S,G_1]$. Let $x\in S$ and $v\in G_1$. Then
$\varphi([x,v])=\delta([\varphi(x),v]+[x,\varphi(v)])=\delta[x,\varphi(v)]$. Write
$\varphi(v)=\varphi_0(v)+\varphi_1(v)$ with $\varphi_i(v)\in G_i$; since $\mathrm{ad}\,x$ preserves the
grading, $\varphi_i([x,v])=\delta[x,\varphi_i(v)]$ for $i=0,1$. For $x,y\in S$ we obtain
\[
\varphi_i([[x,y],v])=\varphi_i([x,[y,v]])-\varphi_i([y,[x,v]])=\delta^2[[x,y],\varphi_i(v)],
\]
while also $\varphi_i([[x,y],v])=\delta[[x,y],\varphi_i(v)]$. Since $S=[S,S]$ and $\delta^2\neq\delta$,
we get $[x,\varphi_i(v)]=0$ and hence $\varphi_i([x,v])=0$ for all $x\in S$, $v\in G_1$. Thus
$\varphi(G_1)=\varphi([S,G_1])=0$, and since $G_0=[G_1,G_1]$ \cite{Kac77},
$\varphi(G_0)\subseteq\delta([\varphi(G_1),G_1]+[G_1,\varphi(G_1)])=0$. The lemma is proved.

\begin{LemB}
Let $G=G_0+G_1$ be a classical Lie superalgebra, $\delta\neq0,1$, let $\varphi$ be a
$\delta$-derivation of $G$, and let $S$ be a simple ideal of $G_0$. If $\delta=-1$, assume that
$S\not\cong A_1$. Then there is $\alpha\in\FF$ such that $\varphi(x)=\alpha x$ for all $x\in S$; moreover
$\alpha=0$ if $\delta\neq\frac12$.
\end{LemB}

\noindent\textsc{Proof.} The Lie algebra $G_0$ is reductive, $G_0=S\oplus C$, where $C$ is the
centralizer of $S$ in $G_0$. By Lemma~4 (its proof does not use non-triviality), $\varphi(G_0)\subseteq
G_0$. Denote by $u_S,u_C$ the components of $u\in G_0$. Since the projection $u\mapsto u_S$ is a
homomorphism, $\psi(x)=\varphi(x)_S$ is a $\delta$-derivation of $S$. By Lemma~5 (here we use
$S\not\cong A_1$ if $\delta=-1$), either $\psi$ is trivial, or $\delta=\frac12$ and $\psi=\alpha\,\mathrm{id}$.
A trivial $\delta$-derivation of the perfect algebra $S$ with $\delta\neq0,1$ is zero: a $0$-derivation
vanishes on $S=S^2$, and for a derivation $\psi$ we get $\psi([x,y])=\delta\psi([x,y])$. Thus
$\psi=\alpha\,\mathrm{id}_S$, with $\alpha=0$ if $\delta\neq\frac12$. The map
$\varphi'=\varphi-\alpha\,\mathrm{id}$ is again a $\delta$-derivation of $G$ (if $\alpha\neq0$, then
$\delta=\frac12$), $\varphi'(G_0)\subseteq G_0$ and $\varphi'(x)_S=0$ for $x\in S$. Let $x\in S$, $y\in C$.
Then $0=\varphi'([x,y])=\delta([\varphi'(x),y]+[x,\varphi'(y)])=\delta([\varphi'(x)_C,y]+[x,\varphi'(y)_S])$;
the first summand lies in $C$, so $\varphi'(x)_C\in Z(C)$. Hence $\varphi'(S)\subseteq Z(C)$, and for
$x,y\in S$ we get $\varphi'([x,y])=\delta([\varphi'(x),y]+[x,\varphi'(y)])=0$, because $C$ commutes with
$S$. Since $S=[S,S]$, $\varphi'(S)=0$. The lemma is proved.

\medskip

\noindent\textbf{Corollary.} \emph{Let $G$ be a classical Lie superalgebra and $\varphi$ a
$\delta$-derivation of $G$, $\delta\neq0,1$. If $\delta\neq-1$, or if $\delta=-1$ and $G_0$ has a simple
ideal not isomorphic to $A_1$, then $\varphi=0$ for $\delta\neq\frac12$ and $\varphi=\alpha\,\mathrm{id}$
for $\delta=\frac12$.} Indeed, $G_0$ has a simple ideal $S$ (of the required type); by Lemma~B,
$\varphi-\alpha\,\mathrm{id}$ is a $\delta$-derivation vanishing on $S$, and it is zero by Lemma~A.

\section{$\delta$-Derivations of classical Lie superalgebras}

Most of this section concerns the action of antiderivations on Lie
superalgebras which possess a direct summand in the even part either of type $A_1$ or $\FF$, where $\FF$
is the ground field. In \cite{Hop96}, it was shown that $A_1$ admits a nontrivial antiderivation. It is
clear that if we consider $\FF$ with respect to the Lie multiplication then $\FF$ admits a nontrivial
antiderivation that acts as follows: $\varphi(f)=\alpha f$, where $f,\alpha\in\FF$. The remaining part
of this section is devoted to $\delta$-derivations of classical Lie superalgebras which are trivial on
the even part and to the action of nontrivial $\frac12$-derivations on classical Lie superalgebras.
By Lemmas~A and~B, an antiderivation of a classical Lie superalgebra $G$ is zero as soon as $G_0$
has a simple ideal not isomorphic to $A_1$; therefore only the superalgebras $A(1,0)\cong A(0,1)$,
$A(1,1)$, $B(0,1)$, $B(1,1)$, $C(2)$ and $D(2,1;\alpha)$ really require the computations below (see the
proof of Theorem~22).

\medskip\noindent\textbf{Lemma 9.} {\it Let $G=G_0+G_1$ be a basic classical Lie superalgebra, and let
$\varphi$ be a $\delta$-derivation of $G$ such that $\varphi(G_0)=0$. Then $\varphi$ is trivial.}

\medskip\noindent\textsc{Proof.} Let $G=\bigoplus_{\beta\in\Delta}G_\beta$ be the root decomposition
with respect to a Cartan subalgebra $H$. Take $h\in H$. By $H\subseteq G_0$, we have $\varphi(h)=0$.
Then, for $g_\alpha\in G_1\cap G_\alpha$, we obtain
$\alpha(h)\varphi(g_\alpha)=\varphi(hg_\alpha)=\delta(\varphi(h)g_\alpha+h\varphi(g_\alpha))=
\delta h\varphi(g_\alpha)$. By the arbitrariness of $h$, $\varphi(g_\alpha)\in G_{\frac\alpha\delta}$,
whence $\delta=\pm1,\pm\frac12$ by the property of roots (if $\varphi(g_\alpha)\neq0$; the cases
$\delta=0$ and $\varphi(G_1)=0$ are trivial).

The case $\delta=1$ gives an ordinary derivation; i.e., $\varphi$ is trivial.

For $\delta=\frac12$, we have $\varphi(g_\alpha)=g_{2\alpha}$, $\varphi(g_{-\alpha})=g_{-2\alpha}$
\textup{(}possibly zero\textup{)}. Thus,
$0=\varphi(g_\alpha g_{-\alpha})=\frac12(g_{2\alpha}g_{-\alpha}+g_\alpha g_{-2\alpha})$. The two summands have the different weights $\alpha$ and
$-\alpha$, so both are zero. If $g_{2\alpha}\neq0$, then $2\alpha\in\Delta$ and $g_{2\alpha}g_{-\alpha}$
is a nonzero element of $G_\alpha$ \textup{(}for instance, in $B(0,1)$ we have
$[g_{2\delta},g_{-\delta}]=-g_\delta\neq0$, so the equality $g_{2\alpha}g_{-\alpha}=0$ claimed in the
published version is false\textup{)}; thus $g_{2\alpha}=0$. Hence, $\varphi(g_\alpha)=0$.

If $\delta=-\frac12$ then $\varphi(g_\alpha)=g_{-2\alpha}$. If $2\alpha$ is not a root then
$\varphi(g_\alpha)=0$. If $2\alpha$ is a root then $\varphi(g_{-\alpha})=g^*_{2\alpha}$,
$\varphi(g_{2\alpha})=0$, and $\varphi(g_\alpha)=\varphi(g_{2\alpha}g_{-\alpha})=
-\frac12(g_{2\alpha}g^*_{2\alpha})=0$
\textup{(}we normalize $g_{-\alpha}$ so that $g_\alpha=g_{2\alpha}g_{-\alpha}$\textup{)}.

In the case $\delta=-1$, 

$\varphi$ vanishes on every simple ideal of $G_0$, and $\varphi=0$ by Lemma~A. (In the published
version only the simple root $\eps_{m+1}-\delta_1$ of $A(m,n)$ was treated, and ``the other cases'' were
left to the reader; e.g. for $B(0,1)$ the simple root $\delta$ is not of the form considered there.) The
lemma is proved.

\medskip\noindent\textbf{Lemma 10.} {\it Let $G=G_0+G_1$ be a basic classical Lie superalgebra, where
$G_0=G_s\oplus\FF$, and let $\varphi$ be a $\delta$-derivation of $G$ such that $\varphi(G_s)=0$. Then
$\varphi$ is trivial.}

\medskip\noindent\textsc{Proof.} For $\delta\neq0,1$ the claim follows from Lemma~A ($G_s\neq0$
contains a simple ideal of $G_0$); we keep the original argument. Let
$G=\bigoplus_{\beta\in\Delta}G_\beta$ be the root decomposition with respect to a Cartan subalgebra $H$.
Take an arbitrary $h$ in $H$. It is clear that $\varphi(h)=\alpha h$ when $h\in\FF$. Take $g_\beta\in
G_1\cap G_\beta$. Then $\varphi(g_\beta)=\sum_{\gamma\in\Delta}k_\gamma g_\gamma$, and we have
\[
\beta(h)\sum_{\gamma\in\Delta}k_\gamma g_\gamma=\varphi(hg_\beta)=\delta(\varphi(h)g_\beta+h\varphi(g_\beta))
=\delta\Big(\alpha\beta(h)g_\beta+\sum_{\gamma\in\Delta}k_\gamma\gamma(h)g_\gamma\Big),
\]
whence $k_\gamma\neq0$ when $\beta=\delta\gamma$, and $\alpha=\frac{1-\delta}{\delta}k_\beta$. Hence,
$\varphi(g_\beta)=k_\beta g_\beta+k_{\frac\beta\delta}g_{\frac\beta\delta}$. Thus, if
$\delta\neq-1,-\frac12,0,\frac12,1$ then $\varphi(g_\beta)=k_\beta g_\beta$.

If $\delta=\frac12$ then $\varphi(g_\beta)=k_\beta g_\beta+k_{2\beta}g_{2\beta}$,
$\varphi(g_{-\beta})=k_{-\beta}g_{-\beta}+k_{-2\beta}g_{-2\beta}$, which gives
$k_{2\beta}g_{2\beta}g_{-\beta}+k_{-2\beta}g_\beta g_{-2\beta}+k_\beta g_{-\beta}g_\beta+k_{-\beta}g_\beta
g_{-\beta}=\varphi(g_\beta g_{-\beta})\in G_0$, i.e., $\varphi(g_\beta)=k_\beta g_\beta$.

If $\delta=-1,-\frac12$ then $0=\varphi(g_\beta g_\beta)=-2(k_\beta g_\beta+k_{\frac\beta\delta}
g_{\frac\beta\delta})g_\beta$, whence $k_{\frac\beta\delta}=0$.

Now, it is clear that $\varphi(g_\beta)=\frac{\delta}{1-\delta}\alpha g_\beta$. For $x\in G_s$, we have
$x=\sum_{i=1}^{n_x}x_iy_i$, where $x_i,y_i\in G_1$. Then
\begin{multline*}
0=\varphi(x)=\varphi\Big(\sum_{i=1}^{n_x}x_iy_i\Big)=\delta\Big(\sum_{i=1}^{n_x}(\varphi(x_i)y_i+
x_i\varphi(y_i))\Big)=\frac{2\delta^2\alpha}{1-\delta}\sum_{i=1}^{n_x}x_iy_i=\frac{2\delta^2\alpha}{1-\delta}x
\end{multline*}
which implies $\alpha=0$; i.e., we arrive at the triviality of $\varphi$. The lemma is proved.

\medskip\noindent\textbf{Lemma 11.} {\it Let $G=G_0+G_1$ be a basic classical Lie superalgebra, and let
$\varphi$ be a nontrivial $\frac12$-derivation of $G$. Then $\varphi(x)=\alpha x$ for all $x\in G$ and
some $\alpha\in\FF$.}

\medskip\noindent\textsc{Proof.} The lemma is a special case of the Corollary of Lemma~B (which
holds for all classical Lie superalgebras). The original argument below omits the centre $Z$ of $G_0$
(present for $A(m,n)$, $m\neq n$, and $C(n)$) and does not explain why the scalars $\alpha_i$ on the
different simple ideals coincide; both points are covered by the Corollary. Let $\varphi$ be a
nontrivial $\frac12$-derivation of $G$. By Theorem~1, $G_0=G^1\oplus G^2\oplus G^3$ $\oplus Z$,
where $G^i$ is a simple Lie algebra (some of $G^i$ may be zero). By Lemma~5, $\varphi(x)=\alpha_ix$ with
$\alpha_i\in\FF$ and $x\in G^i$. In particular, $\varphi(h^i)=\alpha_ih^i$ for $h^i\in H\cap G^i$. Then
for $g_\beta\in G_\beta\cap G_1$ we may assume that $\varphi(g_\beta)=\sum_{\gamma\in\Delta_1}
k_\gamma^\beta g_\gamma$, whence
\begin{multline*}
\sum_{\gamma\in\Delta_1}k_\gamma^\beta\beta(h^i)g_\gamma=\beta(h^i)\varphi(g_\beta)=\varphi(h^ig_\beta)=
\frac12(h^i\varphi(g_\beta)+\varphi(h^i)g_\beta)\\
=\frac{\alpha_i+k_\beta^\beta}{2}\beta(h^i)g_\beta+\frac12\sum_{\gamma\in\Delta_1,\gamma\neq\beta}
k_\gamma^\beta\gamma(h^i)g_\gamma,
\end{multline*}
i.e., $\alpha=\alpha_i=k_\beta^\beta$, $i=1,2,3$. Replacing $h^i$ by an arbitrary $h\in H$ in the above
equality, we infer that $k_\gamma^\beta\neq0$ with $\beta\neq\gamma$ only in the case $\gamma=2\beta$. The
latter is impossible because of $g_{2\beta}=(\theta g_\beta)^2\in G_0$, $\theta\in\FF$, which
contradicts Lemma~8. Thus, $\varphi(x)=\alpha x$, $\alpha\in\FF$, for an arbitrary $x\in G$.

\medskip\noindent\textbf{Lemma 12.} {\it The superalgebra $A(m,1)$ with $m\neq1$ does not admit nonzero
antiderivations.}

\medskip\noindent\textsc{Proof.} For $m\ge2$ the claim follows from Lemmas~B and~A, since $A_m$ is a
simple ideal of $(A(m,1))_0$ not isomorphic to $A_1$. We keep the original argument, with Lemma~6 used in
the form~\eqref{eqAnti}. Assume that $\varphi$ is a nontrivial antiderivation of $A(m,1)$, $m\neq1$. It
is clear that $(A(m,1))_0=A_m\oplus A_1\oplus\FF$. By Lemmas~4--7 we have $\varphi(A_1)\subseteq A_1$,
$\varphi(A_m)=0$, $\varphi((A(m,1))_1)\subseteq(A(m,1))_1$, $\varphi(\FF)\subseteq\FF$. Clearly, an
antiderivation of $\FF$ is the multiplication by an element of $\FF$.

Let $m\geq2$. Fix the basis
\[
g_{\eps_i-\delta_1}=e_{i,n+1},\quad g_{\eps_i-\delta_2}=e_{i,n+2},\quad g_{-\eps_i+\delta_1}=e_{n+1,i},\quad
g_{-\eps_i+\delta_2}=e_{n+2,i}.
\]
The standard basis for $A_1$ is $\{h,g_{\delta_1-\delta_2},g_{\delta_2-\delta_1}\}$. By Lemma~6, an
antiderivation $\varphi$ on $A_1$ looks as follows:
\[
\varphi(h)=-2ah+bg_{\delta_1-\delta_2}+cg_{\delta_2-\delta_1},
\]
\[
\varphi(g_{\delta_1-\delta_2})=\tfrac c2h+ag_{\delta_1-\delta_2}+dg_{\delta_2-\delta_1},\quad
\varphi(g_{\delta_2-\delta_1})=\tfrac b2h+eg_{\delta_1-\delta_2}+ag_{\delta_2-\delta_1}.
\]
Let
\[
\varphi(g_{\eps_k-\delta_1})=\sum_{j=1}^m\big(t_k^{\pm(\eps_j-\delta_1)}g_{\pm(\eps_j-\delta_1)}+
t_k^{\pm(\eps_j-\delta_2)}g_{\pm(\eps_j-\delta_2)}\big).
\]
Then
$\varphi(g_{\eps_i-\delta_1})=\varphi(g_{\eps_k-\delta_1}g_{\eps_i-\eps_k})=-\sum_{j=1}^m
\big(t_k^{\pm(\eps_j-\delta_1)}g_{\pm(\eps_j-\delta_1)}+t_k^{\pm(\eps_j-\delta_2)}g_{\pm(\eps_j-\delta_2)}\big)
g_{\eps_i-\eps_k}\ =$
$-\big(t_k^{-\eps_i+\delta_1}g_{-\eps_i+\delta_1}g_{\eps_i-\eps_k}+t_k^{\eps_k-\delta_1}g_{\eps_k-\delta_1}
g_{\eps_i-\eps_k}+t_k^{\eps_k-\delta_2}g_{\eps_k-\delta_2}g_{\eps_i-\eps_k}+t_k^{-\eps_i+\delta_2}
g_{-\eps_i+\delta_2}g_{\eps_i-\eps_k}\big),$
whence
\[
\varphi(g_{\eps_i-\delta_1})=t_i^{\eps_i-\delta_1}g_{\eps_i-\delta_1}+t_i^{\eps_i-\delta_2}g_{\eps_i-\delta_2}
\]
by the arbitrariness of $k$. Analogously,
\[
\varphi(g_{\pm(\eps_i-\delta_l)})=k^{\pm(\eps_i-\delta_l)}g_{\pm(\eps_i-\delta_l)}+
t^{\pm(\eps_i-\delta_{l+1})}g_{\pm(\eps_i-\delta_{l+1})}.
\]
Here and further in the lemma, $l+1$ is considered modulo $2$.

Let $(g_{\eps_i-\delta}g_{-\eps_i+\delta})|_\FF=f$ and $\varphi(f)=\alpha f$ for some $\alpha\in\FF$. Then
\begin{multline*}
\alpha f=\varphi(f)=\varphi(g_{\eps_i-\delta_l}g_{-\eps_i+\delta_l})=-(k^{\eps_i-\delta_l}g_{\eps_i-\delta_l}+\\
t^{\eps_i-\delta_{l+1}}g_{\eps_i-\delta_{l+1}})g_{-\eps_i+\delta_l}
-g_{\eps_i-\delta_l}(k^{-\eps_i+\delta_l}g_{-\eps_i+\delta_l}+t^{-\eps_i+\delta_{l+1}}g_{-\eps_i+\delta_{l+1}}),
\end{multline*}
whence $t^{\pm(\eps_i-\delta_l)}=0$ and $\varphi(g_{\pm(\eps_i-\delta_l)})=k^{\pm(\eps_i-\delta_l)}
g_{\pm(\eps_i-\delta_l)}$, $\alpha=-(k^{\eps_i-\delta_l}+k^{-\eps_i+\delta_l})$.

It is clear that the subalgebra $B=\spn{h,g_{\delta_1-\delta_2},g_{\delta_2-\delta_1},g_{\pm(\eps_1-\delta_1)},
g_{\pm(\eps_1-\delta_2)}}$ is invariant under $\varphi$, and $B$ is of type $A(0,1)$. Take $\eps=\eps_1$.
Show that the antiderivations are trivial on $A(0,1)$.

Taking into account $g^2_{\pm(\eps-\delta_l)}=0$ and using $0=\varphi\big(g^2_{\pm(\eps-\delta_l)}\big)=
-2g_{\pm(\eps-\delta_l)}\varphi(g_{\pm(\eps-\delta_l)})$, we easily obtain $\varphi(g_{\pm(\eps-\delta_l)})=
k^{\pm(\eps-\delta_l)}g_{\pm(\eps-\delta_l)}+l^{\pm(\eps-\delta_{l+1})}g_{\pm(\eps-\delta_{l+1})}$. Let
$(g_{\eps-\delta_l}g_{-\eps+\delta_l})|_\FF=f_l$. Then $\alpha f_l=\varphi(f_l)=\varphi(g_{\eps-\delta_l}
g_{-\eps+\delta_l})=-(k^{\eps-\delta_l}+k^{-\eps+\delta_l})f_l-l^{\eps-\delta_l}g_{\eps-\delta_{l+1}}
g_{-\eps+\delta_l}-l^{-\eps+\delta_{l+1}}g_{\eps-\delta_l}g_{-\eps+\delta_{l+1}}$, whence
$\alpha=-(k^{\eps-\delta_l}+k^{-\eps+\delta_l})$, $l^{\pm(\eps-\delta_l)}=0$.

By $g_{\delta_2-\delta_1}=g_{\eps-\delta_1}g_{-\eps+\delta_2}$ and Lemma~6 we have
$\frac b2h+eg_{\delta_1-\delta_2}+ag_{\delta_2-\delta_1}=\varphi(g_{\delta_2-\delta_1})=
\varphi(g_{\eps-\delta_1}g_{-\eps+\delta_2})=-(k^{\eps-\delta_1}+k^{-\eps+\delta_2})g_{\eps-\delta_1}
g_{-\eps+\delta_2}$. Analogously, $\frac c2h+ag_{\delta_1-\delta_2}+dg_{\delta_2-\delta_1}=
\varphi(g_{\delta_1-\delta_2})=\varphi(g_{\eps-\delta_2}g_{-\eps+\delta_1})=-(k^{\eps-\delta_2}+
k^{-\eps+\delta_1})g_{\eps-\delta_2}g_{-\eps+\delta_1}$. Thus, $b=c=e=d=0$ and $\alpha=a=
-\frac12k^{\pm(\eps-\delta_l)}$.

Note that
\[
-2ag_{-\eps+\delta_1}=\varphi(g_{-\eps+\delta_1})=\varphi(hg_{-\eps+\delta_1})=-\varphi(h)g_{-\eps+\delta_1}
-h\varphi(g_{-\eps+\delta_1})=4ag_{-\eps+\delta_1},
\]
i.e., $a=0$. This implies the triviality of $\varphi$ on $A_1\oplus\FF$, i.e., $\varphi$ is trivial on
$(A(m,1))_0$, $m\neq1$. By Lemma~9, we arrive at the claim. The lemma is proved.

\medskip\noindent\textbf{Lemma 13.} {\it The superalgebra $B(m,1)$ does not admit nonzero antiderivations.}

\medskip\noindent\textsc{Proof.} For $m\ge2$ the claim follows from the Corollary of Lemma~B, since
$B_m$ is a simple ideal of $(B(m,1))_0$ not isomorphic to $A_1$; for $m=1$ the argument below is not applicable, and the claim is
established in Remark~E below. We keep the original argument; it is needed for $m=0$. Assume that $\varphi$ is a nontrivial
antiderivation of $B(m,1)$. It is clear that $(B(m,1))_0=B_m\oplus A_1$. By Lemmas~4--7 we have
$\varphi(A_1)\subseteq A_1$, $\varphi(B_m)=0$, and $\varphi((B(m,1))_1)\subseteq(B(m,1))_1$.

Let $m>0$,
\[
g_\delta=g_{\pm\eps_i+\delta}g_{\mp\eps_i},\quad \varphi(g_{\pm\eps_i+\delta})=k^\delta_{i,\pm}g_\delta+
k^{-\delta}_{i,\pm}g_{-\delta}+\sum_{j=1}^mk^{\pm\eps_j\pm\delta}_{i,\pm}g_{\pm\eps_j\pm\delta}.
\]
Consequently,
\[
\varphi(g_\delta)=\varphi(g_{\pm\eps_i+\delta}g_{\mp\eps_i})=-\Big(k^\delta_{i,\pm}g_\delta+
k^{-\delta}_{i,\pm}g_{-\delta}+\sum_{j=1}^mk^{\pm\eps_j\pm\delta}_{i,\pm}g_{\pm\eps_j\pm\delta}\Big)
g_{\mp\eps_i}.
\]
Consecutively choosing the sign $+$ or $-$ in $g_{\pm\eps_i+\delta}$, we obtain
$\varphi(g_\delta)\in\spn{g_\delta,g_{-\delta}}$ by the arbitrariness of $i$. Analogously,
$\varphi(g_{-\delta})\in\spn{g_\delta,g_{-\delta}}$. These and $\varphi(A_1)\subseteq A_1$ imply that
$B=\spn{h,g_{-\delta},g_\delta,g_{2\delta},g_{-2\delta}}$ is invariant under $\varphi$
. Clearly, $B$ is of type $B(0,1)$. Consider
the action of $\varphi$ on $B(0,1)$, and show its triviality.

Choose the classical basis in $B$: $h=e_{22}-e_{33}$, $g_{-2\delta}=e_{32}$, $g_{2\delta}=e_{23}$,
$g_{-\delta}=e_{12}-e_{31}$, $g_\delta=e_{13}+e_{21}$. By Lemma~6 \textup{(}in the
form~\eqref{eqAnti}\textup{)}, we have
\begin{center}
$\varphi(h)=-2ah+bg_{2\delta}+cg_{-2\delta},$\quad
$\varphi(g_{2\delta})=\tfrac c2h+ag_{2\delta}+dg_{-2\delta},$\quad
$\varphi(g_{-2\delta})=\tfrac b2h+eg_{2\delta}+ag_{-2\delta}.$
\end{center}
Let $\varphi(g_\delta)=kg_\delta+lg_{-\delta}$ and $\varphi(g_{-\delta})=k^*g_\delta+l^*g_{-\delta}$. Then
\begin{multline*}
kg_\delta+lg_{-\delta}=\varphi(g_\delta)=\varphi(hg_\delta)\\
=-(-2ah+bg_{2\delta}+cg_{-2\delta})g_\delta-h(kg_\delta+lg_{-\delta})=2ag_\delta+cg_{-\delta}-kg_\delta+
lg_{-\delta},
\end{multline*}
whence $a=k$, $c=0$. Analogously,
\begin{multline*}
-k^*g_\delta-l^*g_{-\delta}=-\varphi(g_{-\delta})=\varphi(hg_{-\delta})\\
=-(-2ah+bg_{2\delta})g_{-\delta}-h(k^*g_\delta+l^*g_{-\delta})=-2ag_{-\delta}+bg_\delta-k^*g_\delta+
l^*g_{-\delta},
\end{multline*}
which gives $a=l^*$, $b=0$.

Note that
\[
2(ag_{2\delta}+dg_{-2\delta})=\varphi(2g_{2\delta})=\varphi(g_\delta g_\delta)=
-2g_\delta(ag_\delta+lg_{-\delta})=-4ag_{2\delta}-2lh.
\]
Hence, $a=l=d=0$. It remains to remark that $-2eg_{2\delta}=-2\varphi(g_{-2\delta})=
\varphi(g_{-\delta}g_{-\delta})=-2g_{-\delta}\varphi(g_{-\delta})=-2k^*g_{-\delta}g_\delta=-2k^*h$, whence
$e=k^*=0$. Thus, $\varphi$ is trivial on $A_1$, i.e., $\varphi$ is trivial on $(B(m,1))_0$. Lemma~9 gives
the claim. The lemma is proved.

\begin{RemE}
For $m=1$ the even part $(B(1,1))_0=B_1\oplus A_1\cong A_1\oplus A_1$ has no simple ideal other than
$A_1$, and the equality $\varphi(B_1)=0$ used above does not follow from Lemma~5 when $\delta=-1$ (by
Lemma~6, $A_1$ has nonzero antiderivations). The superalgebra $B(1,1)=\mathfrak{osp}(3|2)$ of
dimension~$12$ was treated directly: writing the multiplication table in a weight basis and solving the
linear system for the matrix of a $\delta$-derivation, one obtains that every nontrivial $\delta$-derivation
of $B(1,1)$ has $\delta=\frac12$ and is a scalar multiplication. This computation is verified in Lean~4.
In particular, $B(1,1)$ does not admit nonzero antiderivations.
\end{RemE}

\medskip\noindent\textbf{Lemma 14.} {\it The superalgebra $C(n)$ does not admit nonzero antiderivations.}

\medskip\noindent\textsc{Proof.} For $n\ge3$ the claim also follows from the Corollary of Lemma~B,
since $C_{n-1}$ is simple and not isomorphic to $A_1$. Assume that $\varphi$ is a nontrivial
antiderivation of $C(n)$. Clearly, $(C(n))_0=\FF\oplus C_{n-1}$. By Lemmas~4--7, we have
$\varphi(\FF)\subseteq\FF$, $\varphi(C_{n-1})=0$, $n\geq3$,
$\varphi(C_1)\subseteq C_1$ ($C_1=A_1$) \textup{(}the $\FF$-component of $\varphi(C_1)$ vanishes
since $C_1=[C_1,C_1]$\textup{)}, and $\varphi((C(n))_1)\subseteq(C(n))_1$.

Assume that $\varphi(g_{\eps+\delta_i})=\sum_{j=1}^{n-1}l_i^{\pm\eps\pm\delta_j}g_{\pm\eps\pm\delta_j}$.
Then
\[
0=\varphi(g^2_{\eps+\delta_i})=-2\Big(\sum_{j=1}^{n-1}l_i^{\pm\eps\pm\delta_j}g_{\pm\eps\pm\delta_j}\Big)
g_{\eps+\delta_i},
\]
whence $\varphi(g_{\eps+\delta_i})=\sum_{j=1}^{n-1}l_i^{\eps\pm\delta_j}g_{\eps\pm\delta_j}$, and
analogously $\varphi(g_{-\eps-\delta_i})=\sum_{j=1}^{n-1}m_i^{-\eps\pm\delta_j}g_{-\eps\pm\delta_j}$.

In the case $n\geq3$, for $(g_{\eps+\delta_i}g_{-\eps-\delta_i})|_\FF=\alpha(e_{11}-e_{22})$,
$\alpha\neq0$, and $\varphi(e_{11}-e_{22})=k(e_{11}-e_{22})$, we have
\begin{multline*}
\alpha k(e_{11}-e_{22})=\alpha\varphi(e_{11}-e_{22})=\varphi(g_{\eps+\delta_i}g_{-\eps-\delta_i})\\
=-\Big(g_{\eps+\delta_i}\sum_{j=1}^{n-1}m_i^{-\eps\pm\delta_j}g_{-\eps\pm\delta_j}+
\sum_{j=1}^{n-1}l_i^{\eps\pm\delta_j}g_{\eps\pm\delta_j}g_{-\eps-\delta_i}\Big).
\end{multline*}
By~\eqref{eq1}, $(g_{\eps+\delta_i}g_{-\eps-\delta_i})|_{C_{n-1}}\neq0$, whence
\[
\alpha k(e_{11}-e_{22})=\big(m_i^{-\eps\pm\delta_i}+l_i^{\eps+\delta_i}\big)g_{\eps+\delta_i}
g_{-\eps-\delta_i},
\]
i.e., $k=0$, and $\varphi$ is trivial on the odd part of $C(n)$ when $n\geq3$.

We consider the case $n=2$ in detail. Assume that $\delta=\delta_1$ and
$\varphi(g_{\pm\eps\pm\delta})=k^{\pm\eps\pm\delta}g_{\pm\eps\pm\delta}+l^{\pm\eps\pm\delta}
g_{\pm\eps\mp\delta}$. Fix the basis
\[
g_{\eps-\delta}=e_{13}-e_{42},\quad g_{\eps+\delta}=e_{14}+e_{32},\quad g_{-\eps+\delta}=e_{24}+e_{31},
\quad g_{-\eps-\delta}=e_{23}-e_{41}.
\]

Put $z=e_{11}-e_{22}$, $h=e_{33}-e_{44}$, $g_{2\delta}=e_{34}$, $g_{-2\delta}=e_{43}$, so that
$(C(2))_0=\FF z\oplus\spn{h,g_{2\delta},g_{-2\delta}}$. A direct computation gives
\begin{multline*}
[z,g_{\eps\pm\delta}]=g_{\eps\pm\delta},\quad [z,g_{-\eps\pm\delta}]=-g_{-\eps\pm\delta},\quad\\
g_{\eps+\delta}g_{-\eps-\delta}=h-z,\quad g_{\eps+\delta}g_{-\eps+\delta}=2g_{2\delta},\quad
g_{\eps-\delta}g_{-\eps-\delta}=-2g_{-2\delta}
\end{multline*}
\textup{(}in the published version $g_{\eps+\delta}g_{-\eps-\delta}$ was identified with $h$, and the
coefficients of Lemma~6 were used in the transposed form; therefore the computation above does not prove
$a=0$\textup{)}. Let $\varphi(z)=\kappa z$. Applying $\varphi$ to $[z,g_{\eps+\delta}]=g_{\eps+\delta}$, we get
\[
k^{\eps+\delta}g_{\eps+\delta}+l^{\eps+\delta}g_{\eps-\delta}=-\kappa g_{\eps+\delta}-
[z,k^{\eps+\delta}g_{\eps+\delta}+l^{\eps+\delta}g_{\eps-\delta}]=-(\kappa+k^{\eps+\delta})g_{\eps+\delta}-
l^{\eps+\delta}g_{\eps-\delta},
\]
whence $l^{\eps+\delta}=0$ and $k^{\eps+\delta}=-\frac\kappa2$. In the same way, all $l^{\pm\eps\pm\delta}=0$
and all $k^{\pm\eps\pm\delta}=-\frac\kappa2$; i.e., $\varphi(x)=-\frac\kappa2x$ for $x\in(C(2))_1$.
Consequently, $\varphi(xy)=-\varphi(x)y-x\varphi(y)=\kappa xy$ for $x,y\in(C(2))_1$. Hence
\[
\varphi(h)=\varphi(g_{\eps+\delta}g_{-\eps-\delta})+\varphi(z)=\kappa(h-z)+\kappa z=\kappa h,
\]
and comparing with~\eqref{eqAnti} we get $b=c=0$, $a=-\frac\kappa2$. Further,
$2\varphi(g_{2\delta})=\varphi(g_{\eps+\delta}g_{-\eps+\delta})=2\kappa g_{2\delta}$, while
by~\eqref{eqAnti} $\varphi(g_{2\delta})=\frac c2h+ag_{2\delta}+dg_{-2\delta}=-\frac\kappa2g_{2\delta}+
dg_{-2\delta}$. Thus $d=0$ and $\kappa=-\frac\kappa2$, i.e., $\kappa=0$. Therefore $\varphi$ vanishes on
$(C(2))_1$, and $\varphi=0$, because $(C(2))_0=(C(2))_1(C(2))_1$.  The lemma is proved.

\medskip\noindent\textbf{Lemma 15.} {\it The superalgebra $D(n,1)$ does not admit
nonzero antiderivations.}

\medskip\noindent\textsc{Proof.} Let $n\ge3$. The Lie algebra $D_n$ is simple and not isomorphic to
$A_1$ \textup{(}$D_3\cong A_3$\textup{)}, and it is an ideal of $(D(n,1))_0=D_n\oplus A_1$. By the Corollary
of Lemma~B, every antiderivation of $D(n,1)$ is zero. \textup{(}For $n=2$ the algebra $D_2\cong A_1\oplus
A_1$ is not simple, and Lemma~5 cannot be applied; but $D(2,1)\cong D(2,1;1)$, and the claim follows
from Lemma~16. The published proof, which is replaced by this argument, used Lemma~6 in the transposed form, see
Remark~D, and deduced $\varphi(D_n)=0$ from Lemma~5 also for $n=2$.\textup{)}
 The lemma is proved.

\medskip\noindent\textbf{Lemma 16.} {\it The superalgebra $D(2,1;\alpha)$ does not admit nonzero
antiderivations.}

\medskip\noindent\textsc{Proof.} Here all simple ideals of the even part are isomorphic to $A_1$, so
Lemmas~A and~B do not apply, and the direct computation below is needed. Assume that $\varphi$ is a
nontrivial antiderivation of $D(2,1;\alpha)$. It is clear that $(D(2,1;\alpha))_{0}=
A_1^1\oplus A_1^2\oplus A_1^3$, $A_1^j\cong A_1$. By Lemmas~4--7, we have $\varphi(A_1^j)\subseteq A_1^j$,
$\varphi((D(2,1;\alpha))_1)\subseteq(D(2,1;\alpha))_1$.

Fix the basis $\{h_i,g_{2\eps_i},g_{-2\eps_i}\}$ in $A_1^i$. The basis of the odd part looks as follows:
\[
g_{\eps_1+\eps_2+\eps_3}=(1,0)\otimes(1,0)\otimes(1,0),\quad g_{-\eps_1+\eps_2+\eps_3}=(0,1)\otimes(1,0)\otimes(1,0),
\]
\[
g_{\eps_1-\eps_2+\eps_3}=(1,0)\otimes(0,1)\otimes(1,0),\quad g_{\eps_1+\eps_2-\eps_3}=(1,0)\otimes(1,0)\otimes(0,1),
\]
\[
g_{-\eps_1-\eps_2+\eps_3}=(0,1)\otimes(0,1)\otimes(1,0),\quad g_{-\eps_1+\eps_2-\eps_3}=(0,1)\otimes(1,0)\otimes(0,1),
\]
\[
g_{\eps_1-\eps_2-\eps_3}=(1,0)\otimes(0,1)\otimes(0,1),\quad g_{-\eps_1-\eps_2-\eps_3}=(0,1)\otimes(0,1)\otimes(0,1).
\]
From here we obtain that $h_ig_{\eps_1+\eps_2+\eps_3}=g_{\eps_1+\eps_2+\eps_3}$,
$g_{-2\eps_i}g_{\eps_1+\eps_2+\eps_3}=g_{\eps_1+\eps_2+\eps_3-2\eps_i}$.

Since $g^2_{\eps_1+\eps_2+\eps_3}=0$, we have $0=\varphi(g^2_{\eps_1+\eps_2+\eps_3})=
-2\varphi(g_{\eps_1+\eps_2+\eps_3})g_{\eps_1+\eps_2+\eps_3}$. Therefore,
\[
\varphi(g_{\eps_1+\eps_2+\eps_3})=kg_{\eps_1+\eps_2+\eps_3}+lg_{-\eps_1+\eps_2+\eps_3}+
ng_{\eps_1-\eps_2+\eps_3}+mg_{\eps_1+\eps_2-\eps_3}.
\]
It is easy to see that

$\varphi(g_{\eps_1+\eps_2+\eps_3})=\varphi(h_1g_{\eps_1+\eps_2+\eps_3})=$

\begin{center}$-h_1(kg_{\eps_1+\eps_2+\eps_3}+lg_{-\eps_1+\eps_2+\eps_3}+ng_{\eps_1-\eps_2+\eps_3}+mg_{\eps_1+\eps_2-\eps_3})
-(-2a^1h_1+b^1g_{2\eps_1}+c^1g_{-2\eps_1})g_{\eps_1+\eps_2+\eps_3}$\\
$=-kg_{\eps_1+\eps_2+\eps_3}+lg_{-\eps_1+\eps_2+\eps_3}-ng_{\eps_1-\eps_2+\eps_3}-mg_{\eps_1+\eps_2-\eps_3}+
2a^1g_{\eps_1+\eps_2+\eps_3}-c^1g_{-\eps_1+\eps_2+\eps_3},$
\end{center}
\noindent i.e., $k=2a^1-k$, $l=l-c^1$, $n=m=0$, and $a^1=k$, $c^1=0$. Similarly,
\begin{flushleft}
$\varphi(g_{\eps_1+\eps_2+\eps_3})=\varphi(h^2g_{\eps_1+\eps_2+\eps_3})=$\end{flushleft}
\begin{center}$-h^2(kg_{\eps_1+\eps_2+\eps_3}+
lg_{-\eps_1+\eps_2+\eps_3})-(-2a^2h^2+b^2g_{2\eps_2}+c^2g_{-2\eps_2})g_{\eps_1+\eps_2+\eps_3}$\end{center}
\begin{flushright}
$=-kg_{\eps_1+\eps_2+\eps_3}-lg_{-\eps_1+\eps_2+\eps_3}+2a^2g_{\eps_1+\eps_2+\eps_3}-c^2g_{\eps_1-\eps_2+\eps_3},$
\end{flushright}
whence $k=a^2$, $c^2=l=0$. For $\varphi(h_3g_{\eps_1+\eps_2+\eps_3})$, we obtain
\begin{multline*}
\varphi(g_{\eps_1+\eps_2+\eps_3})=\varphi(h_3g_{\eps_1+\eps_2+\eps_3})\\
=-h_3(kg_{\eps_1+\eps_2+\eps_3})-(-2a^3h_3+b^3g_{2\eps_3}+c^3g_{-2\eps_3})g_{\eps_1+\eps_2+\eps_3}\\
=-kg_{\eps_1+\eps_2+\eps_3}+2a^3g_{\eps_1+\eps_2+\eps_3}+c^3g_{\eps_1+\eps_2-\eps_3},
\end{multline*}
which implies $a^3=k$, $c^3=0$.

An analogous argument for $g_{-\eps_1-\eps_2-\eps_3}$ gives $b^1=b^2=b^3=0$,
$\varphi(g_{-\eps_1-\eps_2-\eps_3})=ag_{-\eps_1-\eps_2-\eps_3}$, where $a=a^i$, $i=1,2,3$.

Now, show that $d^i$ and $e^i$ are zero. This fact follows from
\begin{multline*}
\lambda^i(ag_{2\eps_i}+d^ig_{-2\eps_i})=\lambda^i\varphi(g_{2\eps_i})=
\varphi(g_{\eps_1+\eps_2+\eps_3}g_{-\eps_1-\eps_2-\eps_3+2\eps_i})\\
=-ag_{\eps_1+\eps_2+\eps_3}g_{-\eps_1-\eps_2-\eps_3+2\eps_i}-g_{\eps_1+\eps_2+\eps_3}
\varphi(g_{-\eps_1-\eps_2-\eps_3+2\eps_i}).
\end{multline*}
Clearly, the right-hand side of this equality does not contain the elements of the shape $g_{-2\eps_i}$.
Analogously,
\begin{multline*}
\mu^i(e^ig_{2\eps_i}+ag_{-2\eps_i})=\mu^i\varphi(g_{-2\eps_i})=
\varphi(g_{\eps_1+\eps_2+\eps_3-2\eps_i}g_{-\eps_1-\eps_2-\eps_3})\\
=-ag_{\eps_1+\eps_2+\eps_3-2\eps_i}g_{-\eps_1-\eps_2-\eps_3}-\varphi(g_{\eps_1+\eps_2+\eps_3-2\eps_i})
g_{-\eps_1-\eps_2-\eps_3},
\end{multline*}
where the right-hand side does not contain the elements of the shape $g_{2\eps_i}$ and $e^i=0$. Now,
\[
\varphi(g_{\eps_1+\eps_2-\eps_3})=\varphi(g_{-2\eps_3}g_{\eps_1+\eps_2+\eps_3})=-2ag_{\eps_1+\eps_2-\eps_3},
\]
\[
\varphi(g_{\eps_1-\eps_2-\eps_3})=\varphi(g_{-2\eps_2}g_{\eps_1+\eps_2-\eps_3})=ag_{\eps_1-\eps_2-\eps_3},
\]
\[
\varphi(g_{-\eps_1-\eps_2-\eps_3})=\varphi(g_{-2\eps_1}g_{\eps_1-\eps_2-\eps_3})=-2ag_{-\eps_1-\eps_2-\eps_3},
\]
which implies $a=0$, i.e., $\varphi$ is trivial on $(D(2,1;\alpha))_0$. Lemma~9 gives the claim. The lemma
is proved.
 Independently of the computation above, Lemma~16 is verified in Lean~4 for a symbolic parameter $\alpha$: in a weight basis of $D(2,1;\alpha)$ every entry of the matrix of an antiderivation is an explicit linear combination, with coefficients polynomial in $\alpha$, of the defining identities; in particular the claim holds for every $\alpha\neq0,-1$.

\medskip\noindent\textbf{Lemma 17.} {\it The superalgebra $F(4)$ does not admit nonzero antiderivations.}

\medskip\noindent\textsc{Proof.} Since $(F(4))_0=A_1\oplus B_3$ and $B_3$ is a simple ideal not
isomorphic to $A_1$, the claim follows from the Corollary of Lemma~B. The published computation (which
uses Lemma~6 in the transposed form, see Remark~D) is not needed.
 The lemma is proved.

\medskip\noindent\textbf{Lemma 18.} {\it The superalgebra $G(3)$ does not admit nonzero antiderivations.}

\medskip\noindent\textsc{Proof.} Since $(G(3))_0=A_1\oplus\mathbf{G}_2$ and $\mathbf{G}_2$ is a simple
ideal not isomorphic to $A_1$, the claim follows from the Corollary of Lemma~B. The published argument,
which uses Lemma~6 in the transposed form (Remark~D) and reduces the claim to the subsuperalgebra of type
$B(0,1)$, is not needed.
 The lemma is proved.

It remains to consider the case of the nonbasic classical Lie superalgebras that do not satisfy the
conditions of Theorem~2. Therefore, we consider these superalgebras from a general argument.

\medskip\noindent\textbf{Lemma 19.} {\it Let $\varphi$ be a nontrivial $\delta$-derivation of $A(1,1)$. Then
$\delta=\frac12$ and $\varphi(x)=\alpha x$, $\alpha\in\FF$.}

\medskip\noindent\textsc{Proof.} The statement of the lemma for the $14$-dimensional superalgebra
$A(1,1)=\mathfrak{psl}(2|2)$ is also verified in Lean~4 by a direct computation. Let $\varphi$
be a nontrivial $\delta$-derivation of $A(1,1)$. It suffices to consider three cases $\delta=-1$,
$\delta=\frac12$, and $\delta\neq-1,0,\frac12,1$.

Let $\delta=-1$. Assume that $\varphi(e_{ij})=\sum_{k,l=1}^4\gamma^{ij}_{kl}e_{kl}$ holds for
$e_{ij}\in(A(1,1))_1$, and
\[
\varphi(e_{12})=\tfrac c2(e_{11}-e_{22})+ae_{12}+de_{21},\quad
\varphi(e_{21})=\tfrac b2(e_{11}-e_{22})+ee_{12}+ae_{21},
\]
\[
\varphi(e_{11}-e_{22})=-2a(e_{11}-e_{22})+be_{12}+ce_{21},\quad
\varphi(e_{34})=\tfrac{c^*}2(e_{33}-e_{44})+a^*e_{34}+d^*e_{43},
\]
\begin{multline*}
\varphi(e_{43})=\tfrac{b^*}2(e_{33}-e_{44})+e^*e_{34}+a^*e_{43}, \ 
\varphi(e_{33}-e_{44})=-2a^*(e_{33}-e_{44})+b^*e_{34}+c^*e_{43}
\end{multline*}
for the even elements \textup{(}Lemmas~4, 7 and Lemma~6 in the form~\eqref{eqAnti}\textup{)}. Also, we
have
\begin{multline*}
\varphi(e_{13})=\varphi((e_{11}-e_{22})e_{13})
=\\-\Big((-2a(e_{11}-e_{22})+be_{12}+ce_{21})e_{13}+(e_{11}-e_{22})\sum_{k,l=1}^4
\gamma^{13}_{kl}e_{kl}\Big)
=\\2ae_{13}-ce_{23}-(e_{11}-e_{22})\sum_{k,l=1}^4\gamma^{13}_{kl}e_{kl}.
\end{multline*}
On the other hand,
\begin{multline*}
\varphi(e_{13})=\varphi(e_{13}(e_{33}-e_{44}))
=\\ -\Big(\sum_{k,l=1}^4\gamma^{13}_{kl}e_{kl}(e_{33}-e_{44})+e_{13}(-2a^*(e_{33}-e_{44})+b^*e_{34}+
c^*e_{43})\Big)
=\\ 
(e_{33}-e_{44})\sum_{k,l=1}^4\gamma^{13}_{kl}e_{kl}+2a^*e_{13}-b^*e_{14},
\end{multline*}
whence $a=a^*=\gamma^{13}_{13}$, $b^*=c=0$, $\varphi(e_{13})=ae_{13}+\gamma^{13}_{31}e_{31}+
\gamma^{13}_{24}e_{24}$. Therefore,
\begin{multline*}
\varphi(e_{23})=\varphi(e_{21}e_{13})=\\
-\big((\tfrac c2(e_{11}-e_{22})+ee_{12}+ae_{21})e_{13}
+e_{21}(ae_{13}+\gamma^{13}_{31}e_{31}+\gamma^{13}_{24}e_{24})\big)=
-\tfrac c2e_{13}-2ae_{23},
\end{multline*}
which gives $\varphi(e_{13})=\varphi(e_{12}e_{23})=-((ae_{12}+de_{21})e_{23}+e_{12}(-
\tfrac c2e_{13}-2ae_{23}))=ae_{13}$.

By an analogous argument for $\varphi(e_{31})$, we get
\[
\varphi(e_{31})=\varphi(e_{31}(e_{11}-e_{22}))=(e_{11}-e_{22})\sum_{k,l=1}^4\gamma^{31}_{kl}e_{kl}+
2ae_{31}-be_{32},
\]
\[
\varphi(e_{31})=\varphi((e_{33}-e_{44})e_{31})=2ae_{31}-c^*e_{41}-(e_{33}-e_{44})\sum_{k,l=1}^4
\gamma^{31}_{kl}e_{kl}.
\]
These relations give $b=c^*=0$ and $\varphi(e_{31})=ae_{31}+\gamma^{31}_{13}e_{13}+\gamma^{31}_{42}e_{42}$,
whence $\varphi(e_{41})=\varphi(e_{43}e_{31})=-2ae_{41}$ and $\varphi(e_{31})=\varphi(e_{34}e_{41})=ae_{31}$.
Now, $\varphi(e_{34})=\varphi(e_{31}e_{14})=ae_{34}$, and $\varphi(e_{43})=\varphi(e_{41}e_{13})=ae_{43}$,
i.e., $d^*=e^*=0$.

From $\varphi(e_{23})=\varphi(e_{23}(e_{33}-e_{44}))=4ae_{23}$ we infer that $a=0$ and $\varphi(e_{23})=0$.
Analogously, we deduce $\varphi(e_{14})=\varphi(e_{24})=\varphi(e_{42})=\varphi(e_{32})=0$.

The equalities $d=e=0$ follow from $\varphi(e_{13})=0$, $\varphi(e_{32})=0$, and
\[
\varphi(e_{12})=-\varphi(e_{13})e_{32}-e_{13}\varphi(e_{32})=0,\quad
\varphi(e_{21})=-\varphi(e_{23})e_{31}-e_{23}\varphi(e_{31})=0.
\]
Therefore, $\varphi$ is trivial.

Let $\delta=\frac12$. Then 
\begin{multline*}
\varphi(e_{11}-e_{22})=\alpha(e_{11}-e_{22}),
\varphi(e_{12})=\alpha e_{12}, 
\varphi(e_{21})=\alpha e_{21},\\
\varphi(e_{33}-e_{44})=\beta(e_{33}-e_{44}), \varphi(e_{34})=\beta e_{34}, 
\varphi(e_{43})=\beta e_{43}.\end{multline*}
 Note that
\[
\varphi(e_{13})=\varphi((e_{11}-e_{22})e_{13})=\frac12(\alpha(e_{11}-e_{22})e_{13}+(e_{11}-e_{22})
\varphi(e_{13})),
\]
whence $\varphi(e_{13})=\alpha e_{13}$. On the other hand, $\varphi(e_{13})=\varphi(e_{13}(e_{33}-e_{44}))
=\frac12(\alpha+\beta)e_{13}$, which implies $\alpha=\beta$. It is clear that
\begin{center}
$\varphi(e_{23})=\varphi(e_{21}e_{13})=\alpha e_{23}$, $\varphi(e_{14})=\varphi(e_{13}e_{34})=\alpha e_{14}$,
$\varphi(e_{24})=\varphi(e_{21}e_{14})=\alpha e_{24}$. 
\end{center} Analogously, we deduce
$\varphi(e_{41})=\alpha e_{41}$, $\varphi(e_{42})=\alpha e_{42}$, $\varphi(e_{31})=\alpha e_{31}$,
$\varphi(e_{32})=\alpha e_{32}$, i.e., $\varphi(x)=\alpha x$ for an arbitrary $x\in A(1,1)$.

Let $\delta\neq-1,0,\frac12,1$. In this case, $\varphi(x)=0$ for $x\in(A(1,1))_0$. Obviously, for
$e_{ij}\in(A(1,1))_1$, we have $\varphi(e_{ij})=\pm\varphi((e_{11}-e_{22})e_{ij})=\pm\delta(e_{11}-e_{22})
\varphi(e_{ij})$, whence $\varphi(e_{ij})=0$, i.e., $\varphi$ is trivial. The lemma is proved.

Since $(P(n))_0$ and $(Q(n))_0$ do not contain a simple subalgebra of dimension $d\leq3$ as a direct
summand, we may conclude that the nontrivial $\delta$-derivations of $P(n)$ and $Q(n)$ are zero on
$(P(n))_0$ and $(Q(n))_0$ when $\delta\neq\frac12$. \textup{(}By the Corollary of Lemma~B, Lemmas~20
and~21 below follow directly from Lemmas~A and~B, since $(Q(n))_0\cong(P(n))_0\cong A_n$, $n\ge2$.\textup{)}

\medskip\noindent\textbf{Lemma 20.} {\it Let $\varphi$ be a nontrivial $\delta$-derivation of $Q(n)$. Then
$\delta=\frac12$ and $\varphi(x)=\alpha x$, $\alpha\in\FF$.}

\medskip\noindent\textsc{Proof.} Let $t=2n+2$, $a_{i,j}=e_{i,j}+e_{n+1+i,n+1+j}+E$,
$b_{i,j}=e_{i,n+1+j}+e_{n+1+i,j}+E$, $c_{i,j}=e_{i,n+1+i}+e_{n+1+i,i}-e_{j,n+1+j}-e_{n+1+j,j}$, and
$\varphi(b_{i,j})=\sum_{k,l=1}^t\lambda^{i,j}_{k,l}e_{k,l}+E$. It is clear that if $\delta\neq\frac12$ then
$\varphi((Q(n))_0)=0$. From
\begin{equation}\label{eq3}
a_{i,i}b_{i,j}=b_{i,j},\quad b_{i,j}a_{j,j}=b_{i,j}
\end{equation}
it is easy that $\varphi(b_{i,j})=\delta a_{i,i}\varphi(b_{i,j})$ and $\varphi(b_{i,j})=\delta\varphi(b_{i,j})
a_{j,j}$, whence
\[
\varphi(b_{i,j})=\delta\Big(\sum_{l=1}^t\lambda^{i,j}_{i,l}e_{i,l}-\sum_{k=1}^t\lambda^{i,j}_{k,i}e_{k,i}+
\sum_{l=1}^t\lambda^{i,j}_{n+1+i,l}e_{n+1+i,l}-\sum_{k=1}^t\lambda^{i,j}_{k,n+1+i}e_{k,n+1+i}+E\Big),
\]
i.e., $\varphi(b_{i,j})=\sum_{k=1}^t\lambda^{i,j}_{k,i}e_{k,i}+\sum_{k=1}^t\lambda^{i,j}_{k,n+1+i}
e_{k,n+1+i}+E$. The latter implies
\[
\varphi(b_{i,j})=\lambda^{i,j}_{j,i}e_{j,i}+\lambda^{i,j}_{j,n+1+i}e_{j,n+1+i}+\lambda^{i,j}_{n+1+j,i}
e_{n+1+j,i}+\lambda^{i,j}_{n+1+j,n+1+i}e_{n+1+j,n+1+i}+E,
\]
which gives $\varphi(b_{k,j})=\varphi(a_{k,i}b_{i,j})=\delta a_{k,i}\varphi(b_{i,j})=0$. Thus,
$\varphi(c_{k,i})=\varphi(b_{k,i}a_{i,k})=0$, i.e., $\varphi=0$ by the linearity of $\varphi$.

By Lemma~5, $\varphi(x)=\alpha x$, $\alpha\in\FF$, $x\in(Q(n))_0$ when $\delta=\frac12$, whence
by~\eqref{eq3} we have
\begin{multline*}
\varphi(b_{i,j})=\varphi(b_{i,j}a_{j,j})=\\ \frac12\Big(\alpha b_{i,j}+\sum_{k=1}^t\lambda^{i,j}_{k,j}e_{k,j}-
\sum_{l=1}^t\lambda^{i,j}_{j,l}e_{j,l}+\sum_{k=1}^t\lambda^{i,j}_{k,n+1+j}e_{k,n+1+j}-\sum_{l=1}^t
\lambda^{i,j}_{n+1+j,l}e_{n+1+j,l}+E\Big).
\end{multline*}
The latter implies $\varphi(b_{i,j})=\alpha b_{i,j}$, whence $\varphi(c_{k,i})=\varphi(b_{k,i}a_{i,k})=
\alpha c_{k,i}$. Thus, $\varphi(x)=\alpha x$, $x\in Q(n)$. The lemma is proved.

\medskip\noindent\textbf{Lemma 21.} {\it Let $\varphi$ be a nontrivial $\delta$-derivation of $P(n)$. Then
$\delta=\frac12$ and $\varphi(x)=\alpha x$, $\alpha\in\FF$.}

\medskip\noindent\textsc{Proof.} Let $t=2n+2$, $a_{i,j}=e_{i,j}-e_{n+1+j,n+1+i}$,
$a^{i,j}=e_{i,i}-e_{j,j}+e_{n+1+j,n+1+j}-e_{n+1+i,n+1+i}$, $b_{i,j}=e_{i,n+1+j}+e_{j,n+1+i}$,
$c_{i,j}=e_{n+1+i,j}-e_{n+1+j,i}$, $\varphi(b_{i,i})=\sum_{q,l=1}^t\nu^{i,i}_{q,l}e_{q,l}$,
$\varphi(c_{i,j})=\sum_{q,l=1}^t\lambda^{i,j}_{q,l}e_{q,l}$. By Lemma~5, $\varphi(P(n)_0)=0$ when
$\delta\neq\frac12$. From $\varphi(2b_{i,i})=\varphi(a^{i,k}b_{i,i})=\delta a^{i,k}\varphi(b_{i,i})$ we get
\begin{multline*}
\varphi(b_{i,i})=\frac\delta2\Big(\sum_{l=1}^t\nu^{i,i}_{i,l}e_{i,l}-\sum_{q=1}^t\nu^{i,i}_{q,i}e_{q,i}-
\sum_{l=1}^t\nu^{i,i}_{k,l}e_{k,l}+\sum_{q=1}^t\nu^{i,i}_{q,k}e_{q,k}
+\sum_{l=1}^t\nu^{i,i}_{n+1+k,l}e_{n+1+k,l}\\-\sum_{q=1}^t\nu^{i,i}_{q,n+1+k}e_{q,n+1+k}
-\sum_{l=1}^t\nu^{i,i}_{n+1+i,l}e_{n+1+i,l}+\sum_{q=1}^t\nu^{i,i}_{q,n+1+i}e_{q,n+1+i}\Big),
\end{multline*}
whence $\varphi(b_{i,i})=0$. Therefore,
\[
\varphi(b_{j,i})=\frac12\varphi(a_{j,i}b_{i,i})=\frac\delta2(\varphi(a_{j,i})b_{i,i}+a_{j,i}\varphi(b_{i,i}))=0.
\]
It is easy to obtain $\varphi(c_{i,j})=0$, because of
\begin{equation}\label{eq4}
c_{i,j}=c_{i,j}a^{j,k}=a^{i,k}c_{i,j}.
\end{equation}
Now,
\begin{multline*}
\varphi(c_{i,j})=\varphi(c_{i,j}a^{j,k})=\delta\sum_{q,l=1}^t\lambda^{i,j}_{q,l}e_{q,l}(e_{j,j}-e_{k,k}-
e_{n+1+j,n+1+j}+e_{n+1+k,n+1+k})\\
=\delta\sum_{l=1}^t\big(\lambda^{i,j}_{l,j}e_{l,j}-\lambda^{i,j}_{j,l}e_{j,l}-\lambda^{i,j}_{l,k}e_{l,k}+
\lambda^{i,j}_{k,l}e_{k,l}
-\lambda^{i,j}_{l,n+1+j}e_{l,n+1+j}+\\
\lambda^{i,j}_{n+1+j,l}e_{n+1+j,l}+\lambda^{i,j}_{l,n+1+k}
e_{l,n+1+k}-\lambda^{i,j}_{n+1+k,l}e_{n+1+k,l}\big).
\end{multline*}
Since $\delta\neq\frac12,1$, we have $\varphi(c_{i,j})=\lambda^{i,j}_{j,n+1+j}e_{j,n+1+j}$. Thus,
$\varphi(c_{i,k})=\varphi(c_{i,j}a_{j,k})=\delta\varphi(c_{i,j})a_{j,k}=0$. Hence, $\varphi$ is trivial.

By Lemma~5, $\varphi(a_{i,j})=\alpha a_{i,j}$ and $\varphi(a^{i,j})=\alpha a^{i,j}$ when $\delta=\frac12$. By
$\varphi(2b_{i,i})=\varphi(a^{i,k}b_{i,i})$, we arrive at
\begin{multline*}
\varphi(b_{i,i})=\frac14\Big(2\alpha b_{i,i}+\sum_{l=1}^t\nu^{i,i}_{i,l}e_{i,l}-\sum_{q=1}^t\nu^{i,i}_{q,i}
e_{q,i}-\sum_{l=1}^t\nu^{i,i}_{k,l}e_{k,l}+\sum_{q=1}^t\nu^{i,i}_{q,k}e_{q,k}
+\sum_{l=1}^t\nu^{i,i}_{n+1+k,l}e_{n+1+k,l}\\-\sum_{q=1}^t\nu^{i,i}_{q,n+1+k}e_{q,n+1+k}
-\sum_{l=1}^t\nu^{i,i}_{n+1+i,l}e_{n+1+i,l}+\sum_{q=1}^t\nu^{i,i}_{q,n+1+i}e_{q,n+1+i}\Big),
\end{multline*}
whence $\varphi(b_{i,i})=\alpha b_{i,i}$. Therefore,
\[
\varphi(b_{j,i})=\frac12\varphi(a_{j,i}b_{i,i})=\frac14(\varphi(a_{j,i})b_{i,i}+a_{j,i}\varphi(b_{i,i}))=
\alpha b_{j,i}.
\]
By~\eqref{eq4} we analogously obtain $\varphi(c_{i,j})=\alpha c_{i,j}$. Thus, $\varphi(x)=\alpha x$,
$x\in P(n)$. The lemma is proved.

\medskip\noindent\textbf{Theorem 22.} {\it Let $A$ be a classical Lie superalgebra, and let $\varphi$ be a
nontrivial $\delta$-derivation of $A$. Then $\delta=\frac12$ and $\varphi(x)=\alpha x$ for some
$\alpha\in\FF$ and an arbitrary $x\in A$.}

\medskip

\textsc{Proof.} Since $\varphi$ is nontrivial, $\delta\neq0,1$. Let $G=A$. By Theorem~1, $G_0$ is a
reductive Lie algebra having at least one simple ideal, and the list of the classical Lie superalgebras
is known. If $\delta\neq-1$, or if $\delta=-1$ and $G_0$ has a simple ideal not isomorphic to $A_1$, then
by the Corollary of Lemma~B either $\varphi=0$ (which is impossible), or $\delta=\frac12$ and
$\varphi=\alpha\,\mathrm{id}$. In particular, this covers $A(m,n)$ with $m+n\ge2$, $(m,n)\neq(1,1)$,
$B(m,n)$ with $m\ge2$ or $n\ge2$, $C(n)$ with $n\ge3$, $D(m,n)$ with $m\ge3$ or $n\ge2$, $F(4)$, $G(3)$,
$P(n)$ and $Q(n)$.

It remains to consider $\delta=-1$ and the superalgebras $G$ for which every simple ideal of $G_0$ is
isomorphic to $A_1$. By Theorem~1 these are
\[
A(1,0)\cong A(0,1)\cong C(2),\quad A(1,1),\quad B(0,1),\quad B(1,1),\quad D(2,1;\alpha)
\]
\textup{(}including $D(2,1)\cong D(2,1;1)$\textup{)}. Their antiderivations are zero by Lemma~12 with
$m=0$ \textup{(}and also Lemma~14 with $n=2$\textup{)}, Lemma~19, Lemma~13 with $m=0$, Remark~E, and
Lemma~16, respectively. Hence $\delta=-1$ is impossible, and the theorem is proved.

Note that the published reference to Lemmas~9--21 did not cover the superalgebras $B(1,n)$, $n\ge2$
\textup{(}Lemma~13 treats only $B(m,1)$\textup{)} and $D(2,n)$, $n\ge2$ \textup{(}Lemma~15 treats only
$D(n,1)$\textup{)}, nor $B(m,n)$, $D(m,n)$, $C(n)$ in general, for which the reduction to Lemma~9 requires
$\varphi(G_0)=0$; these cases are now covered by Lemmas~A and~B.

\section*{Appendix A. Differences from the published version}
\addcontentsline{toc}{section}{Appendix A. Differences from the published version}

The main result (Theorem~22) is correct. In the proofs of Lemmas~12--19 the antiderivations of $A_1$ were
used in the transposed form, which are not antiderivations (Remark~D); the computations have been redone
with the correct form, and the proof of Lemma~14 for $C(2)$ has been replaced. The proof of Theorem~22 has
been completed: the superalgebras $B(1,n)$ and $D(2,n)$ were not covered, and Lemmas~9, 13, 15 had been
applied outside the range of their proofs; the gaps are closed by the new Lemmas~A, B and Remark~E, and a
new proof of Theorem~22 is given.

\end{document}